\documentstyle[10pt]{amsart}

\RequirePackage{amssymb}
\RequirePackage[T1]{fontenc}

\newtheorem{theorem}{Theorem}[section]

\newtheorem{proposition}[theorem]{Proposition}
 
\theoremstyle{definition}

\newtheorem{remark}[theorem]{Remark}

\newcommand{\nc}{\newcommand}
\nc{\Symm}{{\on{Sym}}}

\newcommand{\on}{\operatorname}

 \nc{\cE}{{\cal E}}

\newcommand{\ZZ}{{\mathbb Z}}

\newcommand{\QQ}{{\mathbb Q}}
\newcommand{\CC}{{\mathbb C}}
\newcommand{\RR}{{\mathbb{R}}}

\nc{\SL}{{\mathfrak sl}}
\nc{\gt}{{\mathfrak gt}}
\nc{\grt}{{\mathfrak grt}}
\nc{\gtm}{{\mathfrak gtm}}
\nc{\grtm}{{\mathfrak grtm}}
\nc{\gtmd}{{\mathfrak gtmd}}
\nc{\grtmd}{{\mathfrak grtmd}}
\nc{\pb}{{\mathfrak{pb}}}
\renewcommand{\t}{{\mathfrak{t}}}
\nc{\vt}{{\ov{\t}}}

\nc{\g}{{\mathfrak{g}}}

\nc{\G}{{\mathfrak{G}}}
\nc{\HH}{{\mathfrak H}}

\newcommand{\n}{{\mathfrak{n}}}

\newcommand{\SG}{{\mathfrak{S}}}
\newcommand{\f}{{\mathfrak{f}}}

\nc{\wh}{\widehat}\nc{\wt}{\widetilde}

\newcommand{\ov}{\overline}

\newcommand{\ben}{\begin{enumerate}}
\newcommand{\een}{\end{enumerate}}

\newcommand{\PP}{{\mathbb{P}}}

\newcommand{\TT}{{\bf T}}

\begin{document}

\title[On some relations between generalized associators] 
{On some relations between generalized associators}

\author{Benjamin Enriquez}
\address{IRMA (CNRS), rue Ren\'e Descartes, F-67084 Strasbourg, France}
\email{enriquez@@math.u-strasbg.fr}

\maketitle

\begin{abstract} 
Let $\Phi$ be the Knizhnik-Zamolodchikov associator and $\Psi_N$ be 
its analogue for $N$th 
roots of 1. We prove a hexagon relation for $\Psi_4$. Similarly to the 
Broadhurst (for $\Psi_2$) and Okuda (for $\Psi_4$) duality relations, 
it relies on the "supplementary" (i.e., non-dihedral) symmetries of 
$\CC^\times - \mu_4(\CC)$ (i.e., the octahedron group $\SG_4$). 
We also derive relations between $\Phi$ and $\Psi_2$, which are 
analogues of equations, found by Nakamura and Schneps, satisfied by 
the image of the morphism $\on{Gal}(\bar\QQ/\QQ)\to \widehat{\on{GT}}$.  
\end{abstract}

\tableofcontents

\section*{Introduction}

For $N\geq 1$, let $\Psi_N(A|b[\zeta],\zeta\in\mu_N(\CC))$ be the generalized 
associator defined as the renormalized holonomy from 0 to 1 of 
\begin{equation} \label{eq:N}
(dH)H^{-1} = ({A\over z} + \sum_{\zeta\in \mu_N(\CC)} 
{{b[\zeta]}\over{z-\zeta}})dz,  
\end{equation}
i.e., $\Psi_N = H_1^{-1}H_0$, where $H_0,H_1$ are the solutions of (\ref{eq:N})
on $]0,1[$ 
such that $H_0(z) \simeq z^A$ when $z\to 0^+$, $H_1(z) \simeq (1-z)^{b[1]}$
when $z\to 1^-$ and $A$, $b[\zeta]$ ($\zeta\in\mu_N(\CC)$) are free variables. 
When $N=1$, we set $b[1] = B$, and $\Psi_1(A|b[1]) = \Phi(A,B)$
is the KZ associator. $\Psi_N$ can be viewed as a generating series for the  
values at $N$th roots of unity of multiple polylogarithms. 

In \cite{E}, we found some relations satisfied by $\Psi_N$ 
(the pentagon and octogon relations). As in \cite{Dr}, 
these relations give rise to a torsor and a graded Lie algebra $\grtmd(N)$. 
We have a Lie algebra morphism $\grtmd(N) \to \grt$, where $\grt$ is the 
graded analogue of the Grothendieck-Teichm\"uller Lie algebra (\cite{Dr}). 
The octogon relation is based on the dihedral symmetries of $\CC^\times -
\mu_N(\CC)$. 

However, for special values of $N$, the automorphism group of 
$\CC^\times - \mu_N(\CC)$ is larger than the dihedral group $D_N$. 
These values are $N = 1,2,4$. The resulting supplementary relations are: 
(a) when $N=1$, the duality and hexagon relations (\cite{Dr}); 
(b) for $N=2$, the Broadhurst duality relation (\cite{Br}); 
(c) for $N=4$, the 
Okuda duality relation (\cite{O}) and the relation of Section \ref{Phi4}. 
In cases (a) and (c), the octogon relation is then a consequence of the 
duality and hexagon relations, as the octogon can be cut out in 
two neighboring hexagons. 

In the second part of the paper (Section \ref{Phi:Psi2}), we derive 
relations between $\Phi$ and $\Psi_2$, which are analogues of 
identities in the image of the morphism 
$\on{Gal}(\bar\QQ/\QQ)\to \widehat{\on{GT}}$ found in \cite{NS}. 

One can check that all these relations give rise to subtorsors of the 
torsors studied in \cite{E} and to a Lie subalgebra $\grtmd'(N) \subset 
\grtmd(N)$ for $2|N$ and $4|N$. However, the morphism $\grtmd'(N)\to \grt$ 
is surjective if the final questions of \cite{Dr} are answered in the 
affirmative. 

\section{A hexagon relation for $\Psi_4$} \label{Phi4}

For $N\geq 1$, we denote by $\f_{N+1}$ the Lie algebra with generators 
$A,C,b[\zeta]$, 
$\zeta\in \mu_N(\CC)$, with only relation $A + \sum_{\zeta\in\mu_N(\CC)} 
b[\zeta] + C = 0$, and by $\hat\f_{N+1}$ its degree completion 
(the generators all have degree $1$). Then 
$\Psi_N$ belongs to the group $\on{exp}(\hat\f_{N+1})$.

If $S\subset \PP^1(\CC)$ is finite, let $V:= 
\on{Span}_{\CC}(b[s],s\in S) / \CC(\sum_{s\in S} b[s])$. 
Let $\Gamma\subset \on{PSL}_2(\CC)$ be the group of all the automorphisms 
preserving $S$. Then $\Gamma$ acts on $V$ by $\sigma(b[s]) = b[\sigma(s)]$, 
and the form $\sum_{s\in S} b[s] \on{d}\on{ln}(z-s)$ is $\Gamma$-invariant. 

In particular, if $S = \{0,\infty\} \cup \mu_N(\CC)$, then $\Gamma$ acts on 
$\hat\f_{N+1}$ and on the set of solutions of (\ref{eq:N}).  
In that case, $\Gamma \supset D_N$, where $D_N$
is the dihedral group of order $N$, and $D_N\subsetneq \Gamma$ iff 
$N = 1,2,4$. 

When $N = 4$, $S$ is the octahedron $\{0,\infty\} \cup \mu_4(\CC)$ and 
$\Gamma = \SG_4$. Then $\Gamma$ is presented by generators $s,t$ and 
relations $s^3 = t^4 =
(st)^2 = 1$, where $s : z \mapsto {{1+iz}\over{1-iz}}$ and $t : z\mapsto iz$
(here $i = \sqrt{-1}$). 
The corresponding automorphisms of $\f_5$ are given by 
$$
s : A \mapsto b[1], \; b[1]\mapsto b[i], \; b[i]\mapsto A, \; 
b[-1] \mapsto b[-i], \; b[-i] \mapsto C, \; C\mapsto b[-1], 
$$
and $t: A\mapsto A$, $C\mapsto C$, 
$b[\zeta]\mapsto b[i\zeta]$ (recall that 
$C = -A - \sum_{\zeta\in\mu_4(\CC)} b[\zeta]$). 

Recall that $\Psi_4 = H_1^{-1}H_0$, where 
$H_0(z),H_1(z)$ are the solutions of  
\begin{equation} \label{diff:eq}
(dH)H^{-1} = \big( {A\over z} + \sum_{\zeta\in \mu_4(\CC)}
{{b[\zeta]}\over {z-\zeta}}\big) dz , 
\end{equation}
on $]0,1[$ with behavior $H_0(z) \simeq z^{A}$ when $z\to 0^+$ and 
$H_1(z) \simeq (1-z)^{b[1]}$ when $z\to 1^-$. 

\begin{proposition}
$\Psi_4$ satisfies the following hexagon relation\footnote{If $\n$ is a 
pronilpotent Lie algebra, $x\in \n$ and $a\in\CC \setminus \RR_-$, then 
$a^x := \on{exp}(x\on{ln}(a))$, where $\on{ln}(a)$ is chosen with 
imaginary part in $]-\pi,\pi[$.}  
\begin{equation} \label{hexagon:Psi4}
(2/i)^A s^2(\Psi_4) (2/i)^{b[i]} s(\Psi_4)
(2/i)^{b[1]} \Psi_4 = 1. 
\end{equation}
\end{proposition}

{\em Proof.} Let $H_{0^+},H_{1^-}, H_{1+i0^+},H_{i+0^+},H_{i-i0^+},
H_{i0^+}$ be the solutions of (\ref{diff:eq}) in $\TT := \{z\in\CC | 
\on{Re}(z)\geq 0, \on{Im}(z)\geq 0, |z|\leq 1\} - \{0,1,i\}$ with asymptotic 
behaviors
$H_{0^+}(z) \simeq z^A$ when $z\to 0^+$,  
$H_{1^-}(z) \simeq (1-z)^{b[1]}$ when $z\to 1^-$,  
$H_{1+i0^+}(z) \simeq ({{z-1}\over i})^{b[1]}$ when $z\to 1+i0^+$,  
$H_{i+0^+}(z) \simeq (z-i)^{b[i]}$ when $z\to i+0^+$,  
$H_{i-i0^+}(z) \simeq ({{i-z}\over i})^{b[i]}$ when $z\to i-i0^+$,  
$H_{i0^+}(z) \simeq (z/i)^A$ when $z\to i0^+$.   

Let us still denote by $H_0,H_1$ the prolongations of $H_0,H_1$
to $\TT$. Then we have $H_{0^+}(z) = H_0(z)$, $H_{1^-}(z) = H_1(z)$, 
$H_{1+i0^+}(z) = s(H_0(s^{-1}(z))) 2^{b[1]}$, 
$H_{i+0^+}(z) = s(H_1(s^{-1}(z)))$, 
$H_{i-i0^+}(z) = s^2(H_0(s^{-2}(z))) 2^{b[i]}$, 
$H_{i0^+}(z) = s^2(H_1(s^{-2}(z))) 2^{-A}$.  

Then we have $H_{1^-} = H_{0^+} \Psi_4^{-1}$, 
$H_{1+i0^+} = H_{1^-} i^{b[1]}$, 
$H_{i+0^+} = H_{1+i0^+} 2^{-b[1]} s(\Psi_4^{-1})$, 
$H_{i-i0^+} = H_{i+0^+} i^{b[i]}$, $H_{i0^+} = H_{i-i0^+}
2^{-b[i]} s^2(\Psi_4^{-1}) 2^{-A}$, $H_{0^+} = H_{i0^+} i^A$. 
These equalities imply the hexagon relation. 
\hfill \qed \medskip 

For completeness, we recall that $\Psi_4$ satisfies the Okuda duality relation
(\cite{O}): 
\begin{equation} \label{duality:Psi4}
st(\Psi_4) = 2^{-A} \Psi_4^{-1} 2^{-b[1]},   
\end{equation}
where $st$ is the automorphism of order $2$ given by 
$A \leftrightarrow b[1]$, $b[i]\leftrightarrow b[-i]$, 
$b[-1]\leftrightarrow C$. In \cite{E}, we showed that $\Psi_4$
satisfies the octogon equation 
$$
\Psi_4^{-1} i^{2b[1]} (st^2s^{-1})(\Psi_4) i^C 
(s^{-1}ts^{-1})(\Psi_4^{-1}) i^{2b[i]} t(\Psi_4) i^A = 1,   
$$
but as mentioned in the Introduction, this equation is a consequence of 
(\ref{hexagon:Psi4}) and (\ref{duality:Psi4}). Together with $\Phi$, 
$\Psi_4$ also satisfies a mixed pentagon equation; $\Psi_4$ satisfies a 
group-likeness condition, and the following 
distributivity relations: 
$$
\delta_{42}(\Psi_4) = \Psi_2, \quad \pi_{42}(\Psi_4) = 2^{b[1]}\Psi_2, \quad 
\delta_{41}(\Psi_4) = \Phi, \quad \pi_{42}(\Psi_4) = 4^{b[1]}\Phi, 
$$
where for $d|4$, $\pi_{4d},\delta_{4d} : \f_{4+1} \to \f_{d+1}$ are defined by 
$\pi_{4d}(A) = d'A$, $\pi_{4d}(b[\zeta]) = b[\zeta^{d'}]$, 
$\delta_{4d}(A) = A$, $\delta_{4d}(b[\zeta]) = b[\zeta]$ if 
$\zeta\in \mu_{d}(\CC)$ and 
$=0$ otherwise; here $d' = 4/d$. The generator $b[1]$ of 
$\f_{1+1}$ is denoted by $B$. 

\begin{remark} (duality for $\Psi_2$)
When $N=2$, $S = \{0,\infty,1,-1\}$ is the square and $\Gamma = D_4$. 
Then $\Gamma$ is presented by generators $\rho,\sigma$ and relations 
$\rho^4 = \sigma^2 = (\sigma\rho)^2 = 1$. 
The inclusion $D_4 \subset \SG_4 \subset \on{PSL}_2(\CC)$
is given by $\sigma \mapsto st$, $\rho \mapsto t^2$. 
In \cite{Br}, formula (127), Broadhurst showed the duality 
relation 
\begin{equation} \label{dual:Psi2}
\sigma(\Psi_2) = 2^{-A} \Psi_2^{-1} 2^{-b[1]}. 
\end{equation}
Explicitly, the involutive automorphism $\sigma$ of $\f_3$ is given by 
$\sigma : A \leftrightarrow b[1]$, $b[-1] \leftrightarrow C= -B-b[1]-b[-1]$. 

Since $\sigma\circ\delta_{42} = \delta_{42} \circ (st)$,  
the Okuda duality relation (\ref{duality:Psi4}) for $\Psi_4$, 
together with the distribution relation $\Psi_2 = 
\delta_{42}(\Psi_4)$, implies the Broadhurst
duality relation (\ref{dual:Psi2}) for $\Psi_2$. 
\end{remark}

\section{Relations between $\Phi$ and $\Psi_2$} \label{Phi:Psi2}

In this section, we set $b_0 = b[1]$, $b_1 = b[-1]$. 

\subsection{The relation $\Phi(A,B) = 2^B \Psi_2(A|B,-A-B) 2^A$}

Let $A,B$ be free noncommutative variables. Recall that 
$\Phi(A,B) = G_1^{-1}G_0$, where $G_0,G_1$ are the solutions of 
\begin{equation} \label{KZ}
(dG)G^{-1} = ({A\over u} + {B\over{u-1}})du, 
\end{equation}
such that $G_0(u) \simeq u^A$ as $u\to 0^+$ and $G_1(u) \simeq 
(1-u)^B$ as $u\to 1^-$. 

If $A,b_0,b_1$ are free noncommutative variables, recall that 
$\Psi_2(A|b_0,b_1) = H_1^{-1} H_0$, where $H_0,H_1$ are the 
solutions of 
$$
(dH)H^{-1} = ({A\over z} + {{b_0}\over{z-1}} + {{b_1}\over{z+1}})dz 
$$
such that $H_0(z) \simeq z^A$ as $z\to 0^+$ and 
$H_1(z) \simeq (1-z)^{b_0}$ as $z\to 1^-$. 

Then $\Phi(A,B)\in \on{exp}(\hat\f_2)$ and $\Psi_2(A|b_0,b_1)\in
\on{exp}(\hat\f_3)$, where $\hat\f_2$ (resp., $\hat\f_3$)
is the topologically free Lie algebra
generated by $A,B$ (resp., $A,b_0,b_1$).  

\begin{proposition}
Let $A,B$ be free noncommutative variables. Then
\begin{equation} \label{NS:1}
\Phi(A,B) = 2^B \Psi_2(A|B,-A-B) 2^A.
\end{equation}
\end{proposition}

{\em Proof.} The Broadhurst duality relation (\ref{dual:Psi2}) may be written 
as $\Psi_2(b_0|A,-A-b_0-b_1) = 2_{-A} \Psi_2^{-1}(A|b_0,b_1) 2^{-b_0}$. 
Then if we substitute $b_0 \mapsto B$, $b_1 \mapsto -A - B$, we get 
$\Psi_2(B|A,0) = 
2^{-A}\Psi_2^{-1}(A|B,-A-B)2^{-B}$. The relation then follows from 
the distribution relation $\Psi_2(B|A,0) = \Phi(B,A)$ and from the duality 
relation $\Phi(B,A) = \Phi(A,B)^{-1}$. 

A direct proof is as follows. 
$\Psi_2(A|B,-A-B) = \tilde H_1^{-1} \tilde H_0$, where $\tilde H_0,\tilde H_1$
are the solutions of 
\begin{equation} \label{eq1}
(d\tilde H)\tilde H^{-1} = ({A\over z} + {B\over {z-1}} - {{A+B}\over{z+1}})dz
\end{equation}
such that $\tilde H_0(z) \simeq z^A$ as $z\to 0^+$, $\tilde H_1(z) \simeq 
(1-z)^B$ as $z\to 1^-$. 

Set $u:= 2z/(z+1)$. Then $({A\over z} + {B\over{z-1}} - 
{{A+B}\over{z+1}})dz = ({A\over u} + {B\over {u-1}})du$. 

If follows that $G(u)$ is a solution of (\ref{KZ}), then 
$\tilde H(z) := G(u(z))$ is a solution of (\ref{eq1}). 
In particular, one checks that $\tilde H_0(z) = G_0(u(z)) 2^{-A}$
and $\tilde H_1(z) = G_1(u(z)) 2^B$. Therefore 
$\Psi_2(A|B,-A-B) = \tilde H_1^{-1}\tilde H_0 = 2^{-B} G_1^{-1}G_0 2^{-A}
= 2^{-B}  \Phi(A,B) 2^{-A}$. 
\hfill \qed \medskip 

\begin{remark} Relation (\ref{NS:1}) is the analogue of relation  
$f(\tau_1,\tau_2^2) = \tau_2^{4\rho_2}f(\tau_1^2,\tau_2^2) \tau_1^{2\rho_2}
(\tau_1\tau_2^2)^{-2\rho_2}$ in Theorem 2.2 of \cite{NS}. Here 
$f\in \hat F_2$ is in the image of $\on{Gal}(\bar\QQ/\QQ) \to 
\hat\ZZ \times \hat F_2 \to \hat F_2$ (in particular, $f\in \hat F_2'$, the
commutator subgroup of $\hat F_2$, and $f(x,y)f(y,x)=1$) and 
$\rho_2\in \hat\ZZ$ depends only on $f$ (it is called a Kummer cocycle of 
$f$ in \cite{NS}). This relation takes place in $\hat B_3$, where $B_3$ 
is the braid group with $3$ strands. 

Since $f(x,y)\in \hat F_2'$, it lies in the kernel of the morphism 
$\hat F_2 \to
\ZZ/2\ZZ$, $x\mapsto \bar 1$, $y\mapsto \bar 0$, so there exists a unique 
$h(X|y_0,y_1)\in \hat F_3$ such that $f(x,y) = h(x^2|y,xyx^{-1})$. 
In the same way that the map $\sigma \mapsto f(x,y)$ corresponds to the 
KZ associator $\Phi(A,B)$, the map $\sigma\mapsto h(X|y_0,y_1)$ corresponds to 
$\Psi_2(A|b_0,b_1)$. 

Let $K_3 \subset B_3$ be the pure braid group 
with $3$ strands. It contains the elements $x_{12} = \tau_1^2$, 
$x_{23} = \tau_2^2$ and $x_{13} = \tau_2\tau_1^2 \tau_1^{-1}$; 
$x_{12}x_{13}x_{23} = x_{23}x_{13}x_{12}$ generates the center 
$Z(B_3) \simeq \ZZ$ of $B_3$ and $B_3/Z(B_3) \simeq F_2$ is freely generated
by the classes of $x_{12}$ and $x_{23}$. The relation from \cite{NS}
is then rewritten 
$h(x_{12}|x_{23},x_{13}) = x_{23}^{2\rho_2}f(x_{12},x_{23}) x_{12}^{\rho_2}
(x_{12}x_{13}x_{23})^{-\rho_2}$, which is a relation in $\hat K_3$. 
The image of this relation in $\hat F_2$
is 
$$
f(x,y) = y^{-2\rho_2} h(x|y,(yx)^{-1}) x^{-\rho_2}.
$$ 
(\ref{NS:1}) is an analogue of this relation. 
\end{remark}

\subsection{The relation $\Phi(A,B) = 4^B \Psi_2(A|2B,-2(A+B)) 4^A$}

\begin{proposition}
We have 
\begin{equation} \label{NS:2}
\Phi(A,B) = 4^B \Psi_2(A|2B,-2(A+B)) 4^A
\end{equation} with the above conventions. 
\end{proposition}

{\em Proof.} $\Psi_2(A|2B,-2(A+B)) = \bar H_1^{-1}\bar H_0$, 
where $\bar H_0,\bar H_1$ are the solutions of 
\begin{equation} \label{eq2}
(d\bar H)\bar H^{-1} = ({A\over z} + {{2B}\over {z-1}} 
- {{2(A+B)}\over{z+1}})dz
\end{equation}
with $\bar H_0(z) \simeq z^A$ for $z\to 0^+$, $\bar H_1(z) \simeq (1-z)^{2B}$
when $z\to 1^-$. 

Set $u:= 4z/(z+1)^2$. Then $({A\over z} + {{2B}\over {z-1}} 
- {{2(A+B)}\over{z+1}})dz = ({A\over u} + {B\over {u-1}})du$. 

As above, it follows that if $G(u)$ is a solution of (\ref{KZ}), then 
$\bar H(z) := G(u(z))$ is a solution of (\ref{eq2}).   
Moreover, the expansions $u\simeq 4z$ as $z\to 0$ and $1-u \simeq (1-z)^2/4$
as $z\to 1$ imply that $\bar H_0(z) = G_0(u(z)) 4^{-A}$ and 
$\bar H_1(z) = G_1(u(z)) 4^B$. 
Then $\Phi(A,B) = G_1^{-1}G_0 = 4^B \bar H_1^{-1}\bar H_0 4^A = 
4^B \Psi_2(A|2B,-2(A+B))4^A$. 
\hfill \qed \medskip 

\begin{remark} As before, relation (\ref{NS:2}) is the analogue of 
$f(\tau_1,\tau_2^4) = \tau_2^{8\rho_2} f(\tau_1^2,\tau_2^2) \tau_1^{4\rho_2}
(\tau_1\tau_2^2)^{-4\rho_2}$ in Theorem 2.2 of \cite{NS}.   
\end{remark}

\subsection{The relation $\Psi_2(t_{12} + t_{34}|t_{23},t_{14}) = 2^{Z-t_{23}}
\Phi_{1/2}^{1,23,4} \Phi^{1,2,3} (\Phi^{12,3,4})^{-1}$}

Set
$$
\Phi_{1/2}(A,B) := G_{1/2}^{-1} G_0,
$$ 
where $G_{1/2}$ is the solution 
of (\ref{KZ}) such that $G_{1/2}(1/2) = 1$. Then $G_0(1-z) = \theta(G_1(z))$
and $G_{1/2}(1-z) = \theta(G_{1/2}(z))$, where $\theta\in \on{Aut(\hat\f_2)}$
is the exchange of $A$ and $B$, so $G_1^{-1}G_{1/2} = 
\theta(\Phi_{1/2}(A,B)^{-1}) = \Phi_{1/2}(B,A)^{-1}$. Therefore 
$$
\Phi(A,B) = \Phi_{1/2}(B,A)^{-1}\Phi_{1/2}(A,B).
$$

Recall that $\t_4$ is the Lie algebra with generators $t_{ij}$, $1\leq i\neq
j\leq 4$, with relations $t_{ij} = t_{ji}$ ($i\neq j$), 
$[t_{ij} + t_{ik},t_{jk}] = 0$ ($i,j,k$ different) and $[t_{ij},t_{kl}] = 0$
($i,j,k,l$ different). 

\begin{proposition}
We have 
\begin{equation} \label{eq:Phi}
\Psi_2(t_{12} + t_{34}|t_{23},t_{14}) = 2^{Z-t_{23}}
\Phi_{1/2}^{1,23,4} \Phi^{1,2,3} (\Phi^{12,3,4})^{-1}.
\end{equation} 
Here $Z = \sum_{i<j}t_{ij}$ is a generator of the center of $\t_4$, 
$\Phi^{12,3,4} = \Phi(t^{13} + t^{23},t^{34})$, $\Phi^{1,2,3} =
\Phi(t_{12},t_{23})$ and $\Phi^{1,23,4} = 
\Phi(t_{12} + t_{13},t_{24} + t_{34})$. 
\end{proposition}

As in \cite{NS}, this equation implies the pentagon equation, 
eliminating $\Psi_2$ from it and from the equation obtained using the 
automorphism $t_{ij} \mapsto t_{5-i,5-j}$ of $\t_4$. 

\medskip 

{\em Proof.} The system 
\begin{equation} \label{system}
(d_z G)G^{-1} = ({{t_{12}}\over z} + {{t_{23}}\over{z-w}} 
+ {{t_{24}}\over{z-1}})dz, \quad 
(d_w G)G^{-1} = ({{t_{13}}\over{w}} + {{t_{23}}\over{w-z}} 
+ {{t_{34}}\over{w-1}})dw
\end{equation}
is known to be compatible. The pentagon identity is established by 
considering the pentagon of asymptotic zones $((0z)w)1$, 
$(0(zw))1$, $0((zw)1)$, $0(z(w1))$, $(0z)(w1)$ in the domain 
$\{(z,w)|0<z<w<1\}$ (\cite{Dr}). We will cut this pentagon in two 
quadrangles, one of which is $((0z)w)1$, $(0(zw))1$, $(z+w = 1, z \to 1/2^-)$, 
$(z+w = 1, z\to 0^+)$. 

Let $G_{((0z)w)1}$ be the solution 
which is $\simeq z^{t_{12}}w^{t_{13} + t_{23}}$ in the zone 
$((0z)w)1$ (i.e., $w\to 0$, $z/w\to 0$). 

Let $G_{(0(zw))1}$ be the solution which is $\simeq (w-z)^{t_{23}}
w^{t_{12} + t_{13}}$ in the zone $(0(zw))1$ (i.e., $w\to 0$, $z/w\to 1$). 

Let $G_{z+w = 1, z \to 1/2^-}$ be the solution which is $\simeq 
(1-2z)^{t_{23}}$ when $z\to 1/2^-$ and $z+w=1$. 

Let $G_{z+w=1,z\to 0^+}$ be the solution which is 
$\simeq z^{t_{12} + t_{34}}$ when $z\to 0^+$, $z+w = 1$. 
Then $G_{z+w=1,z\to 0^+} = G_{(0z)(w1)}$, where 
$G_{(0z)(w1)}$ is the solution which is $\simeq z^{t_{12}}
(1-w)^{t_{34}}$ when $z\to 0^+$, $w\to 1^-$. 

We know from \cite{Dr} that $G_{(0(zw))1}^{-1} G_{((0z)w)1} 
= \Phi^{1,2,3}$, $G_{((0z)w)1}^{-1} G_{(0z)(w1)} = (\Phi^{12,3,4})^{-1}$. 

Let us compute $G_{(0(zw))1}^{-1}G_{z+w=1,z\to 1/2^-}$. 
If $G$ is a solution of (\ref{system}), 
set $\Gamma(z) := [(w-t)^{-t_{23}} G(z,w)]_{w=z}$. Then 
$$
(d\Gamma)\Gamma^{-1} = ({{t_{12}+t_{13}}\over{z}} 
+ {{t_{24}+t_{34}}\over{z-1}})dz, 
$$
so $[(w-z)^{-t_{23}} G_{(0(zw))1}(z,w)]_{w=z} = G_0(z)^{1,23,4}$, 
$[(w-z)^{-t_{23}} G_{z+w=1,z\to 1/2^-}(z,w)]_{w=z} = G_{1/2}(z)^{1,23,4}$. 
Therefore $G_{(0(zw))1}^{-1}G_{z+w=1,z\to 1/2^-} = (\Phi_{1/2}^{-1})^{1,23,4}$. 

We now compute $G_{z+w=1,z\to 1/2^-}^{-1}G_{z+w=1,z\to 0^+}$. 
If $G$ is a solution of (\ref{system}), set $\Lambda(z) := G(z,1-z)$. 
Then 
$$
(d\Lambda)\Lambda^{-1} = ({{t_{12}+t_{34}}\over z} + {{t_{23}}\over{z-1/2}} 
+ {{t_{13}+t_{24}}\over{z-1}})dz. 
$$
If we set $u:= z/(1-z)$, this equation is 
$$
(d\Lambda)\Lambda^{-1} = ({{t_{12}+t_{34}}\over u} 
+ {{t_{23}}\over{u-1}} + {{t_{14}-Z}\over {u+1}}). 
$$
Then the expansions $u\simeq z$ then $z\to 0$, $1-u \simeq 2(1-2z)$
when $z\to 1/2$ give
$$
\Lambda_{z+w,z\to 0^+}(z) = \hat H_0(u(z)), \quad 
\Lambda_{z+w,z\to 1/2^-}(z) = \hat H_1(u(z)) 2^{-t_{23}},  
$$
where $\hat H_0,\hat H_1$ are the analogues of $H_0,H_1$
for $(A,b_0,b_1) = (t_{12}+t_{34},t_{23},t_{14}-z)$. 
Then 
\begin{align*}
& G_{z+w=1,z\to 1/2^-}^{-1} G_{z+w=1,z\to 0^+} = 
\Lambda_{z+w=1,z\to 1/2^-}^{-1} \Lambda_{z+w=1,z\to 0^+} = 
2^{t_{23}}\hat H_1^{-1} \hat H_0 
\\ & = 2^{t_{23}}\Psi(t_{12}+t_{34}|t_{23},t_{14}-Z)
= 2^{t_{23}-Z}\Psi(t_{12}+t_{34}|t_{23},t_{14}). 
\end{align*}

Then 
$$
G^{-1}_{((0(zw))1} G_{(0z)(w1)} =
(G^{-1}_{((0(zw))1} G_{((0z)w)1})
(G^{-1}_{((0z)w)1} G_{(0z)(w1)}) =
\Phi^{1,2,3} (\Phi^{-1})^{12,3,4}, 
$$
on the other hand 
\begin{align*}
& G^{-1}_{((0(zw))1} G_{(0z)(w1)} =
(G^{-1}_{((0(zw))1} G_{z+w=1,z\to 1/2^-})
(G^{-1}_{z+w=1,z\to 1/2^-} G_{(0z)(w1)}) \\ & =
(\Phi_{1/2}^{-1})^{1,23,4} 2^{t_{23}-Z}\Psi_2(t_{12}+t_{34}|t_{23},t_{14}) .  
\end{align*}
The result follows from the comparison of these equalities. 
\hfill \qed \medskip

\begin{remark} Relation (\ref{eq:Phi}) is the analogue of the relation 
(III') in \cite{NS} 
$$
f(\tau_1\tau_3,\tau_2^2) =
g(x_{45},x_{51})f(x_{12},x_{23})f(x_{34},x_{45}),
$$ 
where $(\lambda,f)\in \hat\ZZ\times \hat F_2$ lies in the image of 
$\on{Gal}(\bar\QQ/\QQ)$ and $g(x,y)\in\hat F_2$
is defined by $g(y,x)^{-1}g(x,y) = f(x,y)$. This relation takes place in 
the uncolored mapping class group $\hat \Gamma_{0,[5]}$ of a surface of 
genus $0$ with $5$ marked points; here $\tau_i,x_{ij}$ are the images of the 
standard elements of the braid group with 5 strands $B_5$ under the morphism 
$B_5 \to \Gamma_{0,[5]}$ (recall that $\tau_i\tau_{i+1}\tau_i = 
\tau_i\tau_{i+1}\tau_i$, $\tau_i\tau_j = \tau_j\tau_i$ if $|i-j| \geq 2$, 
$x_{ij} = \tau_{j-1}...\tau_{i+1} \tau_i^2 (\tau_{j-1}...\tau_{i+1})^{-1}$
if $i<j$, $x_{ji} = x_{ij}$). Indeed, $\tau_2^2 = x_{34}$, 
$(\tau_1\tau_3)^2 = x_{12}x_{34}$ and $(\tau_1\tau_3)\tau_2^2 
(\tau_1\tau_3)^{-1} = x_{24}^{-1}x_{14}x_{24}$, so that  
$f(\tau_1\tau_3,\tau_2^2) = h(x_{12}x_{34}|x_{23},x_{24}^{-1}x_{14}x_{24})$, 
hence (III') is rewritten as 
$$
h(x_{12}x_{34}|x_{23},x_{24}^{-1}x_{14}x_{24})
= g(x_{45},x_{51})f(x_{12},x_{23})f(x_{34},x_{45}) 
$$
which is now an equation in the colored mapping class group $\Gamma_{0,5}$. 
Now $x_{45} = x_{12}x_{13}x_{23}$ and $x_{51} = x_{23}x_{24}x_{34}$, 
so $g(x_{45},x_{51}) = x_{23}^\alpha g(x_{12}x_{13},x_{24}x_{34})$, 
where $\alpha\in\hat\ZZ$, and $f(x_{34},x_{45}) = f(x_{34},x_{13}x_{23})$
(no $x_{12}^\alpha$ comes out since $f\in \hat F_2'$), and using $f(x,y) = 
f(y,x)^{-1}$ the equation is rewritten as 
\begin{equation} \label{final:pent}
h(x_{12}x_{34}|x_{23},x_{24}^{-1}x_{14}x_{24})
= x_{23}^\alpha g(x_{12}x_{13},x_{24}x_{34})
f(x_{12},x_{23}) f(x_{13}x_{23},x_{34})^{-1}. 
\end{equation}
This is an equality in the image of the morphism $\hat K_4 \to 
\hat\Gamma_{0,5}$, where $K_4$ is the pure braid group with $4$ strands. 
This image is the quotient of $\hat K_4$ by its center, 
generated by $x_{12}x_{13}x_{14}x_{23}x_{24}x_{34}$. So the image of 
(\ref{eq:Phi}) in $\on{exp}(\hat\t_4/\CC Z)$ is the analogue of 
(\ref{final:pent}). 

The sense in which the relations of this section are analogous to 
relations in \cite{NS} can be precised as in \cite{Dr}. 
\end{remark}

\end{document}